\documentclass[12pt]{article}
\usepackage{fullpage}
\usepackage{amsfonts, amsmath, amssymb,latexsym,epsfig}
\usepackage{amsthm}
\usepackage{tikz}
\usetikzlibrary{graphs}
\usetikzlibrary{graphs.standard}
\usepackage{relsize}
\usepackage{verbatim}
\usepackage{enumerate}
\usepackage[skip=10pt plus1pt, indent=20pt]{parskip}

\usepackage{hyperref}

\newtheorem{thm}{Theorem}[section]

\newtheorem{lem}[thm]{Lemma}

\newtheorem{prop}[thm]{Proposition}

\DeclareMathOperator{\PG}{PG}

\newcommand{\F}{\mathbb F}

\newcommand{\hs}{\hspace{1em}}

\author{
Jozefien D'haeseleer \thanks{Department of Mathematics: Analysis, Logic and Discrete Mathematics, Ghent University, 9000 Ghent, Belgium. E-mail: {\tt jozefien.dhaeseleer@ugent.be}.}
\and
Vladislav Taranchuk \thanks{Department of Mathematics: Analysis, Logic and Discrete Mathematics, Ghent University, 9000 Ghent, Belgium. E-mail: {\tt vlad.taranchuk@ugent.be}.}
}
\title{On the Chromatic Number of Grassmann Graphs}

\begin{document}

\maketitle

\begin{abstract}

In this paper we study the chromatic number of the Grassmann graphs $J_q(n, m)$. We show that $\binom{n-m+1}{1}_q \leq \chi(J_q(n, m)) \leq \binom{n}{1}_q$, which is analogous to the best-known bounds for the chromatic number of the Johnson graphs $J(n, m)$. When $m = 2$, determining $\chi(J_q(n, 2))$ is equivalent to determining the smallest number of partial line parallelisms that one can partition the lines of $\PG(n-1, q)$ into. We survey known results about line parallelisms and their implications for $\chi(J_q(n, 2))$. Finally, we prove that when $q$ is any power of two, and $n$ is any even integer, then $\chi(J_q(n, 2)) < 2\binom{n-1}{1}_q$. 
   
\end{abstract}

\section{Introduction}

Let $G$ be a graph. The chromatic number of $G$, denoted $\chi(G)$, is the smallest integer $c$ for which $G$ can be partitioned into $c$ independent sets. If $\Delta$ is the maximum degree of $G$, then it is well-known that $\chi(G) \leq \Delta + 1$. In this paper, we will study the chromatic number of a well-known family of graphs giving new bounds. Throughout the paper, $q$ will always represent some prime power. The finite field of order $q$ will be denoted by $\F_q$ and the vector space of dimension $n$ over $\F_q$ will be denoted $\F_q^n$. We use the standard notation $\PG(n, q)$ to denote the projective space of dimension $n$ over the field $\F_q$. We remind the reader that any space of projective dimension $k$ corresponds to a vector space of dimension $k+1$.\\ The gaussian binomial coefficient $\binom{n}{m}_q$ is defined as 
$$
\binom{n}{m}_q = \frac{(q^n - 1)(q^{n-1}-1)\cdots(q^{n-m+1}-1)}{(q-1)(q^2 - 1) \cdots (q^m - 1)}.
$$

 It is a well-known fact that $\binom{n}{m}_q$ represents the total number of distinct $m$-dimensional subspaces of $\F_q^n$ (or equivalently, the total number of distinct $(m-1)$-spaces in $PG(n-1, q)$). For more preliminary background on projective geometries, we refer the reader to \cite{hirshfeld}.

Let $n \geq m$. The Johnson graph $J(n, m)$ is the graph whose vertex set is given by all $m$-subsets of $[n] := \{ 1, 2, \dots, n\}$ and two vertices are adjacent if their intersection has cardinality $m-1$. It is now folklore that 
$$
n - m + 1 \leq \chi(J(n, m)) \leq n.
$$
The lower bound comes from the fact that any $(m-1)$-subset is contained in $n - m+1$ $m$-subsets, each of which must receive distinct colors in any valid coloring. On the other hand, originally shown by Grahame and Sloane \cite{GS}, a valid coloring can be given by assigning to an $m$-subset $\{a_1, a_2, \dots, a_m \}$ the color $a_1 + a_2 + \cdots +a_m \pmod{n}$. Any two distinct $m$-subsets intersecting in an $(m-1)$-subset clearly must have distinct colors, and so this does indeed yield a valid coloring. Here we describe pairs $(n, m)$ for which $\chi(J(n, m))$ is known, or is known to have better bounds than the general bound stated above. 

\begin{enumerate}
    \item For all even $n$, $\chi(J(n, 2)) = n-1$, see \cite{GS}.
    \item For all odd $n$, $\chi(J(n, 2)) = n$, see \cite{GS}.
    \item For $n \leq 6$, $\chi(J(n, 3)) = n$ and $J(7, 3) = 6$.
    \item For $n > 7$, and $n \equiv 1 \pmod{6}$ or $n \equiv 3 \pmod{6}$, $\chi(J(n, 3)) = n-2$, see \cite{Lu1983, Lu1984, Teirlinck1991}.
    \item For $n > 7$, and $n \equiv 0 \pmod{6}$ or $n \equiv 2 \pmod{6}$, $\chi(J(n, 3)) = n-1$, see \cite{Lu1983, Lu1984, Teirlinck1991}.
    \item For $n \equiv 4 \pmod{6}$ and $n \equiv 5\pmod{6}$ it is known that for some infinite sequences of $n$ satisfying these conditions, $J(n, 3) = n$, as shown by Etzion \cite{Etzion1992, Etzion1992P}. It is conjectured that this is correct for all $n > 7$ and $n \equiv 4 \pmod{6}$ and $n \equiv5\pmod{6}$.
\end{enumerate}
Some further improvements to $\chi(J(n, m))$ in certain special cases were also proven by Etzion and Bitan in \cite{Etzion1996}.

Similar to the Johnson graph, the Kneser graph $K(n, m)$ has all  $m$-subsets of $[n]$ as its vertices and two vertices are adjacent if their intersection is empty.
The chromatic number of the Kneser graph $K(n,m)$ was determined by Lov\'asz in 1978 \cite{LovaszKneser}. For $n<2m$ the Kneser graph $K(n,m)$ has no edges, and hence, its chromatic number equals $1$. For $n\geq 2m$ Lov\'asz proved that $\chi(K(n,m))=n-2m+2$. A coloring $C$ with this number of colors can be defined as follows. For
$i \in \{1,\dots, n\}$ let $C_i$ be the class that contains the (ordered) $m$-subsets whose first element is $i$. Then $C_1, C_2, \dots , C_{n-2m+1}$ and $C_{n-2m+2} \cup \dots \cup C_n$ form a partition of the $m$-subsets into $n-2m+2$ classes. 

The $q$-Kneser graph $K_q(n,m)$ is the $q$-analogue of the Kneser graph, in which the vertices are the vector $m$-spaces in the vector space $\mathbb{F}^n_q$ (or equivalently, the $(m-1)$-spaces in the projective space $\PG(n-1,q)$). Two vertices are adjacent if the corresponding $m$-spaces meet trivially.

%In the vector space setting, the subspaces of dimension $1$ and $2$ correspond to the $0$ and $1$ dimensional projective subspaces, and are called the points and lines in this space.

Here again, for $n<2m$ the graph $K_q(n,m)$ has no edges, so its chromatic number is $1$. 
For $m\geq 3$ and $n\geq 2m+1$ (except for $n=2m+1$ and $q=2$) the following theorem is known.
\begin{thm}[{Blokhuis et. al. \cite[Theorem~1.5]{HM2010}}]
Suppose \(m \geq 3\) and either \(q \geq 3\) and \(n \geq 2m+1\) or \(q=2\) and \(n \geq 2m+2\), then the chromatic number of the $q$-Kneser graph is \(\chi(K_q(n,m)) = \binom{n-m+1}{1}_q\). Moreover, each color class of a minimum coloring is a set of $m$-spaces through a fixed point, and the points determining a color are the points of an \((n-m)\)-dimensional subspace of \(\mathbb{F}_q^n\).
\end{thm}
For the remainder values of $m$ and $n$ the following results are known:
\begin{itemize}
    \item Let $m=2$, and $n>4$ then \(\chi(K_q(n,2)) = \binom{n-1}{1}_q = \binom{n-m+1}{1}_q\), as shown by Chowdhury, Godsil, and Royle \cite[Theorem~8.1]{Chowdhury2006}.
    \item Let $m=2$ and $n=4$ then \(\chi(K_q(4,2)) = q^2+ q\), see \cite[Theorem~6.3]{Chowdhury2006}.
    \item  Let \(n=2m\) and \(q\geq 5\), then \(\chi(K_q(2m,m)= q^{m} +q^{m-1}\), as shown by Ihringer \cite{Ihringer2019}.
\end{itemize}

Denote by $J_q(n, m)$ the Grassmann graph where the set of vertices is given by the set of all vector $m$-spaces of $\F_q^n$ and two vertices are adjacent if they intersect in an $(m-1)$-space. Hence, for $m = 2$, we recover the line incidence graph of $\PG(n-1, q)$. Observe that if it is possible to color $J_q(n,2)$ with $\binom{n-1}{1}_q$ colors, then this is equivalent to a partition of the lines of $\PG(n-1, q)$ into spreads. Such a partition of the lines of $\PG(n-1, q)$ is called a line parallelism. It turns out that the existence of line parallelisms has been a topic of research for many years in the realm of finite geometry.

It is known that line spreads can only exist in $\PG(n-1, q)$ when $n-1$ is odd. Indeed, for each odd $n-1$ and each prime power $q$, spreads are known to exist. This of course gives a necessary condition for the existence of line parallelisms, namely that $n-1$ must be odd as well. However, it is still an open question to determine whether line parallelisms exist for each pair $(n-1, q)$. To date, there have been no non-existence results for any $q$ when $n-1$ is odd. However, existence has been determined in some general cases which we outline below. For more information about spreads and parallelisms in $\PG(n-1, q)$, see the survey by Johnson \cite{Johnson2010}.

\begin{enumerate}
    \item Line parallelisms exist in $\PG(n-1, q)$ when $q = 2$ and for each odd $n-1$ as shown by Baker \cite{Baker76}. This implies $\chi(J_q(n, 2)) = \binom{n-1}{1}_q$ for these $n$ and $q$.
    \item Line parallelisms exist in $\PG(n-1, q)$ when $q$ is any prime power and $n = 2^k$ for any positive integer $k \geq 2$ as shown by Beutelspacher \cite{Beutelspacher74}. This implies $\chi(J_q(n, 2)) = \binom{n-1}{1}_q$ for these $n$ and $q$.
    \item Line parallelisms exist in $\PG(n-1, q)$ for each odd $n-1$ when $q = 3, 4, 8, 16$ as shown by Feng, Xu \cite{Xu2023}. This implies $\chi(J_q(n, 2)) = \binom{n-1}{1}_q$ for these $n$ and $q$.
\end{enumerate}
Although spreads do not exist in $\PG(n-1, q)$ when $n-1$ is even, the chromatic number of $J_2(n, 2)$ is known. In \cite{Meszka2013}, Meszka determined that $\chi(J_2(n, 2)) = \binom{n-1}{1}_2 + 3$ when $n-1$ is even.

The Grassmann graph $J_q(n, m)$ is a strongly regular graph with valency $d = q\binom{m}{1}_q\binom{n-m}{1}_q$. To our knowledge, no general bound which improves upon the trivial $d + 1$ bound has been given in the literature for the chromatic number of Grassmann graphs $J_q(n, m)$. Here, we prove two results. The first bounds are general bounds which are analogous to the folklore general bounds for Johnson graphs.

\begin{thm}\label{MainThm1}
Let $q$ be a prime power and $m, n$ be integers with $m < n$. Then 
$$
\binom{n-m+1}{1}_q\leq \chi(J_q(n, m)) \leq \binom{n}{1}_q.
$$
\end{thm}

A quick comparison with the valency of $J_q(n, m)$ shows that for fixed $q$ and $m$ and growing $n$, this bound improves upon the trivial upper bound by a multiplicative factor of $(q-1)/q$. We prove these bounds in Section \ref{sec:bound1}. \\

Next, we give an order of magnitude improvement to the general upper bound for $J_q(n, 2)$ when $q = 2^e$ and $n$ is even.

\begin{thm}\label{MainThm2}
Let $e$ be a positive integer, $q = 2^e$, and $n$ be any even integer. Then 
$$
\binom{n-1}{1}_q \leq \chi(J_q(n, 2)) < 2\binom{n-1}{1}_q.
$$
\end{thm}

This bound is proven in Section \ref{sec:bound2}.

\section{A general bound for $\chi(J_q(n, m))$}\label{sec:bound1}

In this section, we will prove Theorem \ref{MainThm1}. The result is obtained by making some key observations about a special determinant. The following proposition is a known fact, but we give a proof for completeness.

\begin{prop}\label{DetL}
    Let $q$ be a prime power and $m < n$ be a positive integer. Let $x_1, x_2, \dots, x_m \in \F_{q^n}$. The determinant 
    $$
    \left|  
    \begin{matrix}
        x_1 & x_2 & \dots & x_m \\
    x_1^q & x_2^q & \dots & x_m^q \\ 
    \vdots & \vdots& \ddots & \vdots \\
    x_1^{q^{m-1}} & x_2^{q^{m-1}} & \cdots & x_m^{q^{m-1}}
    \end{matrix}
    \right| = 0
    $$
    if and only if $x_1, x_2, \dots, x_m$ are linearly dependent as vectors over $\F_q$.
\end{prop}

\begin{proof}
    Denote by $M$ the matrix in the proposition.
    Suppose that the vectors $x_1, \dots, x_m$ are linearly dependent over $\F_q$. Then there exist scalars $a_1, \dots, a_m \in \F_q$, not all zero for which.
    $$
    a_1x_1 +  a_2x_2 + \cdots a_mx_m = 0.
    $$
    This implies that 
    $$
    0 =  (a_1x_1 +  a_2x_2 + \cdots a_mx_m)^{q^i}= a_1x_1^{q^i} +  a_2x_2^{q^i} + \cdots a_mx_m^{q^i} = 0, \hspace{1em} \forall i \in \mathbb{N}.
    $$
    Thus, the matrix $M$, as matrix over $\F_{q^n}$ has a non-trivial null space, and therefore has determinant 0.

    On the other hand, suppose that the determinant of the above matrix is 0. Consider $M^T$, and suppose the $v = (v_1, \dots, v_m)$ is in its null space. This implies that the function
    $$
    f(x) = v_mx^{q^{m-1}} + \cdots + v_1x
    $$
    has $x_1, \dots, x_m$ as roots. Since $f(x)$ is a linearized polynomial, it means that its set of roots contains the full span of $x_1, \dots, x_m$ as vectors over $\F_q$. Since it's degree is $q^{m-1}$, then it can have at most $q^{m-1}$ roots, therefore, the span of $x_1, \dots, x_m$  over $\F_q$ is at most $(m-1)$-dimensional implying a linear dependence amongst $x_1, \dots, x_m$ over $\F_q$.
\end{proof}

Recall the following properties of determinants.

\begin{prop}\label{DetProp}
    Let $A$ be a square matrix:
    \begin{enumerate}
        \item If $B$ is obtained from $A$ by multiplying a row (or column) by $\alpha$, then $\det(B) = \alpha\det(A)$.
        \item If $B$ is obtained from $A$ by adding to row $j$ (or column $j$) any linear combination of the remaining rows (or columns), then $\det(A) = \det(B)$.
        \item Let $A, B, C$ be matrices such that all rows (or columns) of $A, B, C$ are identical except for row $j$ (or column $j$). If row $j$  (or column $j$) of $C$ is equal to the sums of rows $j$ (or columns $j$) in $A$ and $B$, then 
        $$
        \det(C) = \det(A) + \det(B)
        $$
    \end{enumerate}
\end{prop}

We are now ready to give a coloring of $J_q(n, m)$ which yields the desired bound. 
\begin{proof}[Proof of Theorem \ref{MainThm1}]

Identify the vector space $\F_q^n$ with the field $\F_{q^n}$ by fixing any choice of basis for $\F_{q^n}$ over $\F_q$. Our set of colors will be the cosets of $\F_{q^n}^*/\F_q^*$.

    Let $S$ be an $m$-dimensional subspace of $\F_{q^n}$. Let $\{x_1, \dots, x_m \}$ be a basis for $S$. Define
   
    $$
    L(x_1, x_2, \dots, x_m) =
    \left[
    \begin{matrix}
        x_1 & x_2 & \dots & x_m \\
    x_1^q & x_2^q & \dots & x_m^q \\ 
    \vdots & \vdots& \ddots & \vdots \\
    x_1^{q^{m-1}} & x_2^{q^{m-1}} & \cdots & x_m^{q^{m-1}}
    \end{matrix}
    \right].
    $$
    Color the subspace $S$ with coset $C$ if $\det(L(x_1, x_2, \dots, x_m)) \in C$. We comment that Proposition \ref{DetL} implies that $\det(L(x_1, \dots, x_m)) \neq 0$ whenever $x_1, \dots, x_m$ are linearly independent, and thus this determinant does necessarily belong to some coset of $\F_{q^n}^*/\F_q^*$. The following claim implies that this is well-defined, that is, the coset to which the $\det(L(x_1, x_2, \dots, x_m))$ belongs to is independent of the choice of basis for $S$.

    \noindent \textbf{Claim}: Let $x_1, \dots, x_m$ and $y_1, \dots, y_m$ be two distinct bases for $S$.
    Then $\det(L(x_1, \dots, x_m) )= \alpha \det(L(y_1, \dots, y_m))$, for some $\alpha \in \F_q^*$.

    Note that we may obtain the matrix $L(y_1, \dots, y_m)$ from the matrix $L(x_1, \dots, x_m)$ via elementary column operations since their columns span the same space. Proposition \ref{DetProp} implies that such operations can only change the determinant by a scalar multiple of $\F_q$. Thus the coloring we suggest is indeed well-defined. 

    Now we show the validity of the coloring. Suppose that $S_1$ and $S_2$ are two distinct $m$-dimensional subspaces intersecting in an $(m-1)$-dimensional subspace $T$. Then note we may give a basis for $S_1$ and for $S_2$ which only differ by one element. Let $\{ x_1, \dots, x_{m-1}, x_m \}$ be a basis for $S_1$ and $\{ x_1, \dots, x_{m-1}, x_m^* \}$ be a basis for $S_2$. Since these two subspace are adjacent in $J_q(n, m)$, we must demonstrate they have distinct colors. Suppose not, and so 
    $$
    \det(L(x_1, \dots, x_{m-1}, x_m) = \alpha \det( L(x_1, \dots, x_{m-1}, x_m^*) )= \det(L(x_1, \dots, x_{m-1}, \alpha x_m^*)).
    $$
    Note thas by Proposition \ref{DetProp} (3) and then (1), we have
    $$
    \det(L(x_1, \dots, x_{m-1}, x_m - \alpha x_m^*)) = \det(L(x_1, \dots, x_{m-1}, x_m)) - \alpha \det(L(x_1, \dots, x_{m-1},  x_m^*)) = 0.
    $$
    This implies that the columns of the matrix on the left are linearly dependent. In particular, there exists some  $a_1, \dots, a_{m-1} \in \F_q$ such that 
    $$
    x_m - \alpha x_m^* = a_1x_1 + \cdots + a_{m-1}x_{m-1}.
    $$
    Thus, $x_m^* \in S_1$, a contradiction, and so $S_1$ and $S_2$ must indeed receive distinct colors. There are $\binom{n}{1}_q$ cosets of $\F_{q^n}^*$, and so we have demonstrated a valid coloring of $J_q(n, m)$ with $\binom{n}{1}_q$ colors. 

    The lower bound follows from the fact that the number of $m$-dimensional subspaces in $\F_q^n$ which contain a given $(m-1)$-dimensional subspace is $\binom{n-m+1}{1}_q$. The collection of these $m$-dimensional subspaces, pairwise intersect in $J_q(n, m)$, and so form a clique of size $\binom{n-m+1}{1}_q$. Consequently, each vertex in this clique should belong to a different color class, and so $\binom{n-m+1}{1}_q \leq \chi(J_q(n, m))$.
    \end{proof}

\section{Coloring $J_q(n, 2)$ when $q = 2^e$ and $n$ even}\label{sec:bound2}
In this section, we give a coloring of $J_q(n, 2)$ for all even $n$ and $q = 2^e$ using no more than
$$
2\binom{n-1}{1}_q
$$
colors. While our coloring works for all even $n$ and $q = 2^e$, we remark that for a few values of $n$ and $q$, known results on line parallelisms imply better results than the one our coloring yields. The key idea of our coloring is demonstrating a homomorphism from $J_q(n, 2)$ into the 2-Kneser graph $K_2((n-1)e, e)$, whose chromatic number is known. This in turn implies that 
$$
\chi(J_q(n, 2)) \leq \chi \left(K_2((n-1)e, e)\right) = \binom{ne-2e+1}{1}_2  < 2\binom{n-1}{1}_q
$$

Recall that the $q$-Kneser graph $K_q(n, m)$ is the graph whose vertices are $m$-dimensional subspaces of $\F_q^n$ and two distinct vertices are adjacent if their corresponding subspaces intersect trivially. The function which will be used to demonstrate the homomorphism has was introduced by Hawtin \cite{DH}, where he proved that $\PG(n, q)$ admits $(q-1)$-fold packings, when $q = 2^{e}$, and $n$ is any odd integer. 

Let $n = 2k + 2$ and let $q = 2^e$. Consider the space $V =  \F_{q^{2k + 1}} \times \F_q \cong \F_{q}^{2k+2}$, with vectors denoted by $\overline{x} = (x, x_1) \in \F_{q^{2k + 1}} \times \F_q$. This isomorphism can be established by fixing any choice of basis for $\F_{{q}^{2k+1}}$ when thought of as a vector space over $\F_q$. 

Define a bivariate function on $V$ as follows
$$
g(\overline{x}, \overline{y}) = (x_1y + y_1x)^{q+1} + xy^q + x^qy. 
$$
We note that $g(\overline{x}, \overline{y}) \in \F_{q^{2k+1}}$.

\begin{lem}\label{Lem2}
\noindent Let $\overline{x}, \overline{y} \in V$ be linearly independent (over $\F_q$) vectors, and denote their span by $S$. For any choice of $\overline{z}, \overline{w} \in S$, we have
$$
g(\overline{z}, \overline{w}) \in \{ \alpha^2(x_1y + y_1x)^{q+1} + \alpha(xy^q + x^qy) : \alpha \in \F_q\} = E.
$$
\end{lem}
\begin{proof}
    We may write $\overline{z} = a_1\overline{x} + a_2\overline{y}$ and $\overline{w} = b_1\overline{x} + b_2 \overline{y}$ with $a_1, a_2, b_1, b_2 \in \mathbb{F}_q$. We write
\begin{align*}
 &g(\overline{z}, \overline{w}) = g(a_1\overline{x} + a_2\overline{y}, b_1\overline{x} + b_2 \overline{y}) \\
 =& [(a_1x_1 + a_2y_1)(b_1x + b_2y) + (b_1x_1 + b_2y_1)(a_1x + a_2y)]^{q+1} \\ &+ (a_1x + a_2y)(b_1x + b_2y)^q + (a_1x + a_2y)^q(b_1x + b_2y) \\
 =& (a_1b_2 + a_2b_1)^2(x_1y + y_1x)^{q+1} + (a_1b_2 + a_2b_1)(xy^q + x^qy).
\end{align*}
Setting $\alpha = a_1b_2 +a_2b_1$ finishes the proof.
\end{proof}

\begin{lem}
    Let $\overline{x}, \overline{y} \in V$ be two linearly independent vectors. Then the set
    $$
    E = \{ \alpha^2(x_1y + y_1x)^{q+1} + \alpha(xy^q + x^qy) : \alpha \in \F_q\}
    $$
    is closed under addition. Furthermore, $|E| = q = 2^e.$
\end{lem}

\begin{proof}
    Observe that given any $\alpha, \beta \in \F_q$
\begin{align*}
    &\alpha^2(x_1y + y_1x)^{q+1} + \alpha(xy^q + x^qy) + \beta^2(x_1y + y_1x)^{q+1} + \beta(xy^q + x^qy) \\
    =& (\alpha + \beta)^2(x_1y + y_1x)^{q+1} + (\alpha + \beta)(xy^q + x^qy)
\end{align*}
Thus $E$ is closed under addition and may indeed be thought of as a subspace of $\F_2^{(2k+1)e}$. 

To prove the second claim, we just need to verfiy that for each $\alpha$, the terms in $E$ are distinct. Since $E$ is closed under addition, we just need to verify that the zero vector is only obtained from $\alpha = 0$ in the definition of $E$. 

Suppose not, and therefore, there exists an $\alpha \in \F_q^*$ such that 
\begin{equation}\label{Eq1}
    \alpha^2(x_1y + y_1x)^{q+1} + \alpha(xy^q + x^qy) = 0.
\end{equation}

If $x_1 = y_1 = 0$, then $xy^q + x^qy = 0$ only has the solution $y = \alpha x$. But we assumed that $\overline{x}$ and $\overline{y}$ are linearly independent, so this cannot happen. 

Without loss of generality, suppose $x_1 \neq 0$. Observe that
\begin{align}\label{eq3}
    x_1^{-1}\left[x\left(x_1 y + y_1x\right)^q + x^q\left(x_1y + y_1 x\right)\right] = xy^q + x^qy
\end{align}
Then (\ref{Eq1}) may be re-written as
\begin{equation}\label{Eq1.1}
(x_1y + y_1x)^{q+1} = \left(\frac{x}{\alpha x_1} \right)(x_1y + y_1x)^q + \left(\frac{x}{\alpha x_1} \right)^q(x_1y + y_1x) 
\end{equation}
We observe that the equation 
\begin{equation}%\label{Eq1.2}
z^{q+1} = wz^q +w^qz
\end{equation}
can be re-written as 
$$
z^{q+1} = z^{q+1} + w^{q+ 1} + (z + w)^{q+1} \iff w^{q+1} = (z+w)^{q+1}.
$$
When $q = 2^e$, $q + 1$ is relatively prime with $q^{2k + 1} - 1$, and so $p(x) = x^{q+1}$ is a permutation polynomial over $\F_{q^{2k+1}}$. Thus the equation only has a solution when $z = 0$. This implies that in (\ref{Eq1.1}), we must have $x_1y + y_1x = 0$ again implying a linear dependence between $\overline{x}$ and $\overline{y}$, a contradiction.
\end{proof}

This demonstrates that $g(\overline{x}, \overline{y})$ maps any  2-dimensional subspaces of $V$ to an $e$-dimensional subspace of $\F_2^{(2k+1)e}$. By abuse of notation, given any  2-dimensional subspace $S$, we denote it's associated $e$-dimensional subspace in $\F_2^{(2k+1)e}$ by $g(S)$. We are left to show our final claim.

\begin{prop}\label{LIprop}
    Let $\overline{x}, \overline{y}, \overline{z} \in V$ be non-zero vectors such that the subspaces $S_1 = \text{span}(\overline{x}, \overline{y})$ and $S_2 = \text{span}(\overline{x}, \overline{z})$ are distinct. Then  $g(S_1) \cap g(S_2) = \{ 0\}$.
\end{prop}

\begin{proof}

    Let $\overline{x}$ be a fixed non-zero vector. The number of distinct 2-dimensional subspaces of $V \cong \F_q^{2k + 2}$ which contain $\overline{x}$ is $(q^{2k + 1} - 1)/(q-1)$. Thus, for a fixed $\overline{x}$, and for each 2-dimensional subspace $S$ of $V$ through $\overline{x}$ we obtain a corresponding $e$-dimensional output subspace $E \subset \F_{2}^{(2k+1)e}$.

    Note that a set $S$ of distinct pairwise trivially intersecting $e$-dimensional subspaces in $\F_2^{(2k+1)e}$ contains at most $(q^{2k + 1} - 1)/(q-1)$ elements, and equality occurs if and only if $S$ is a spread and therefore covers each vector in $\F_2^{(2k+1)e}$. This implies that in order to show that the outputs of $g(S)$ at a fixed point $\overline{x}$ intersect trivially, we just need to show that each vector $v \in \F_2^{(2k+1)e}$ belongs to $g(S)$ for some $S$ which contains $\overline{x}$. Thus it is sufficient to demonstrate that for any $\overline{x} \neq 0$, the function $g(\overline{x}, \overline{y})$, as a function of $\overline{y}$, is onto.

    Recall that 
    $$
    g(\overline{x}, \overline{y}) = (x_1 y + y_1 x)^{q+1} + xy^q + x^qy
    $$
    
    \noindent \textbf{Case 1}: Suppose $x_1 \neq 0$, then observe that by linearity and equation \eqref{eq3}, we may re-write $g$ as follows:
    $$
    g(\overline{x}, \overline{y}) = (x_1 y + y_1 x)^{q+1} + x_1^{-1}[x(x_1y + y_1 x)^q + x^q(x_1y + y_1 x)].
    $$
    Performing a change of variables $z = x_1 y + y_1 x$ and $\lambda = x_1^{-1}$, we obtain a new function 
    $$
    f(z) = z^{q+1} + \lambda (xz^q + x^q z).
    $$
    To verify the claim, we need to show that $f(z)$ is permutation of $\F_{q^{2k + 1}}$ for each $x \in \F_{q^{2k + 1}}$ and $\lambda \in \F_q^*$. Note
    $$
     f(z) = z^{q+1} + \lambda (xz^q + x^q z) =  z^{q+1} + (\lambda x)z^q + (\lambda x)^q z = (z + \lambda x)^{q + 1} + (\lambda x)^{q+1}.
    $$
    We remark again that $p(x) = x^{q+1}$ is a permutation polynomial over $\F_q$, and thus $f(z)$ is indeed a permutation for any fixed $\lambda \in \F_q^*$ and $x \in \F_{q^{2k + 1}}$.

    \noindent \textbf{Case 2}: Suppose that $x_1 = 0$ and therefore $x \neq 0$, then 
    $$
    g(\overline{x}, \overline{y}) = (y_1 x)^{q+1} + xy^q + x^qy.
    $$
    Note that for a fixed $x \neq 0$ $xy^q + x^qy$ can be seen to be a linear operator on $\F_{q^{2k + 1}}$ over the field $\F_q$. It quickly follows that this operator has a 1-dimensional kernel. Hence, the image of $xy^q + x^qy$ has dimension $2k$, and so, size $q^{2k}$.  We will show that for any two distinct choices of $y_1 \in \F_q$, call them $a$ and $b$, the functions
    $$
    h_a(y) = a^2x^{q+1} + xy^q + x^qy, \hs  h_b(y) = b^2x^{q+1} + xy^q + x^qy
    $$
    have disjoint images. Suppose for a contradiction that they share a common element, so that 
    $$
    a^2x^{q+1} + xy^q + x^qy = b^2x^{q+1} + xz^q + x^qz \implies (a+b)^2x^{q+1} = x(y+z)^{q} + x^q(y+z).
    $$
    Recall that in the proof of Lemma \ref{Lem2}, we showed that the equation
    $$
    z^{q+1} = wz^q + w^qz
    $$
    only has the solution $z = 0$. Since the image of $wz^q + w^qz$ is an $\F_q$-linear subspace of $\F_{q^{2k+1}}$, then it is also true that 
    $$
    cz^{q+1} = wz^q + w^qz
    $$
    only has the solution $z = 0$ for any $c \in \F_q^*$. Consequently,
    $$
    (a+b)^2x^{q+1} = x(y+z)^{q} + x^q(y+z)
    $$
    only has solutions when $x = 0$ or $a = b$, both of which are contradictions to our assumptions.
    Thus, in case 2, we also obtain that $g(\overline{x}, \overline{y})$ is onto.

    This complete the proof of the proposition.
\end{proof}

\noindent Consequently, when $q = 2^e$, the function $g$ is a homomorphism from $J_q(2k + 2, 2)$ into $K_2((2k+1)e, e)$ which implies that 
$$
\chi(J_q(2k + 2, 2)) \leq \chi(K_2((2k+1)e, e)) = \binom{2ke + 1}{1}_2 < 2 \binom{2k+1}{1}_q.
$$

\section{Concluding Remarks}

Our upper bound for the chromatic number $J_q(n, 2)$, when $q = 2^e$ and $n$ even, improves upon the known bounds except when $q = 2, 4, 8, 16$, or when $n$ is itself a power of 2.

When $q = 2^e$ and $n$ is even, the homomorphism we have from $J_q(n, 2)$ into $K_2((n-1)e, e)$ suggests a route for further improving the upper bound. In particular, we observe that the image of the homomorphism maps specifically into the subset of $e$-spaces which have the form 
$$
\{\alpha^2 (x_1y + y_1x)^{q+1} + \alpha(xy^q + x^qy) : \alpha \in \F_q\}
$$
for any $\overline{x}, \overline{y} \in \F_{q^{n-1}}\times \F_q$. 
Thus, if one can determine the chromatic number of the subgraph of $K_2((n-1)e, e)$ induced by the $e$-spaces of the form above, then one would improve the the upper bound for $J_q(n, 2)$ as well. 

However, having run some computations for $q = 4$ and $n = 4$, it seems likely that the chromatic number of this induced subgraph is strictly larger than the lower bound for the chromatic number of the corresponding Grassmann graphs. Therefore, pursuing this route would not ultimately yield an upper bound matching the lower bound, which is known to be the true answer in many cases.

\bibliographystyle{abbrv}
\bibliography{sample.bib}

\end{document}